\documentclass{amsart}
\usepackage{amsmath,fullpage,amsthm}
\usepackage{amsfonts,color}
\usepackage{latexsym}
\usepackage{url}
\usepackage{amscd}
\usepackage{easybmat}

\newtheorem{thm}{Theorem}[section]

\newtheorem{prop}[thm]{Proposition}
\newtheorem{lemma}[thm]{Lemma}
\newtheorem{cor}[thm]{Corollary}

\newtheorem*{question*}{Question}
\newtheorem{question}[thm]{Question}

\newtheorem{defin}[thm]{Definition}
\newtheorem{observ}[thm]{Observation}
\newtheorem{example}[thm]{Example}
\newtheorem{remark}[thm]{Remark}

\newcommand{\restrictto}[2]{\left. #1 \right|_{#2}}
\newcommand{\ddtat}{\restrictto{\frac{d}{dt}}{t=0}}
\newcommand{\J}{ \begin{BMAT}{cc}{cc}\ 0&1\\ -1&0\end{BMAT}} 

\title{Ricci Yang-Mills solitons on nilpotent Lie groups}
\author{Michael Jablonski}
\address[Michael Jablonski]{ University of Oklahoma}
\email{mjablonski@math.ou.edu}

\author{Andrea Young}
\address[Andrea Young]{ University of Arizona}
\email{ayoung@math.arizona.edu}

\long\def\symbolfootnote[#1]#2{\begingroup%
\def\thefootnote{\fnsymbol{footnote}}\footnote[#1]{#2}\endgroup}

\begin{document}
\maketitle

\begin{abstract}
The purpose of this paper is to introduce the Ricci Yang-Mills soliton equations on nilpotent Lie groups.  In the 2-step nilpotent setting, we show that these equations are strictly weaker than the Ricci soliton equations.  Using techniques from Geometric Invariant Theory, we develop a procedure to build many different kinds of Ricci Yang-Mills solitons.  We finish this note by producing examples of Lie groups that do not admit Ricci soliton metrics but that do admit Ricci Yang-Mills soliton metrics.
\end{abstract}

\section{Introduction}

The goal of this paper is to introduce Ricci Yang-Mills solitons on nilpotent Lie groups.  In this setting, Ricci Yang-Mills solitons are weaker than Ricci soliton metrics in a sense to be made precise below.  We provide some examples of manifolds known not to admit Ricci solitons that do admit Ricci Yang-Mills solitons.

To study this problem, we rephrase our questions in the language of moment maps for a particular representation of $GL_n\mathbb R$.  A similar implementation of Geometric Invariant Theory has been carried out by Lauret, Eberlein, Jablonski, et.~al, in the study of Ricci solitons on nilpotent Lie groups (see, for example, \cite{La},  \cite{Eber07}, \cite{Jab:Thesis}).  Moreover, our approach to the study of Ricci Yang-Mills solitons answers a technical question asked by Eberlein concerning moment maps in the 2-step nilpotent setting (cf. Question \ref{question: distinguished pts of GLqR action}).\\


The Ricci flow is a differential equation on the space of Riemannian metrics on $M$, $\mathfrak{Met}$.  In this space, the only fixed points of this equation are the Ricci-flat metrics, whereas general Einstein metrics are fixed points of the volume normalized Ricci flow.  However, if one works in the space of $\mathfrak{Met}$/$\mathfrak{Diff}$, where $\mathfrak{Diff}$ is the group of diffeomorphisms on $M$, then one allows for a new family of fixed points, namely the metrics that flow by scaling and diffeomorphism; i.e. $g(t)=\sigma(t)\phi(t)^*g_0$, where $\phi(t):M\to M$ is a one parameter family  of diffeomorphisms.    These are the Ricci soliton metrics.  One can show that Ricci soliton metrics satisfy the following equation:
\begin{equation}
\label{soliton}
Rc+\mathcal{L}_Xg+\frac{\epsilon}{2}g=0,
\end{equation}
where $X$ is the vector field generating the diffeomorphisms, and $\epsilon=-1, 0, 1$ corresponds to shrinking, steady, and expanding solitons, respectively.  If $X$ is the gradient of some function, i.e. $X=\nabla f$, then a solution to Eq.~\ref{soliton} is said to be a gradient Ricci soliton.

Ricci soliton structures on Lie groups were first discovered by Baird and Danielo [BD] and independently by Lott [Lo].   Baird and Danielo discovered the first known examples of nongradient soliton structures on the Lie groups Nil and Sol [BD].  They studied semiconformal maps from 3-manifolds to Riemann surfaces and  described completely the soliton structures on all the 3-dimensional geometries.  More generally, one knows quite a bit of information about left-invariant Ricci solitons on Lie groups.  If $G$ is a semi-simple group, then any Ricci soliton metric must be Einstein, and  all compact semi-simple groups admit Einstein metrics.  Within the set of solvable groups, it is known that there exist solvable $G$ which admit non-Einstein Ricci soliton metrics. There also exist solvable groups that do admit Einstein metrics.  Further restricting to nilpotent Lie groups, it is known that there are no Einstein metrics on $G$, but there can be Ricci soliton metrics on $G$ (see \cite{La} for more details).

However, there are some spaces that are known not to admit even Ricci solitons.  Consider a Lie group $G$ with left-invariant metric.   Recently there have been many new families of nilpotent Lie groups constructed which cannot admit left-invariant Ricci soliton metrics \cite{Jab:Moduli}.  To better understand these spaces, it would be useful to have an even weaker notion of `best' metric.  It is our intention to show that a suitable notion could be that of a Ricci Yang-Mills soliton.

The Ricci Yang-Mills flow was defined independently in \cite{St} and \cite{Yo}.   These equations are motivated by the study of Ricci flow on principal bundles and can be written as a modified Ricci flow coupled to the Yang-Mills heat flow.   Long-time behavior of the flow has been studied in \cite{St1}, and stability properties have been considered in \cite{Yo1}.

Gradient Ricci Yang-Mills solitons have been studied in  \cite{St}, \cite{St1}, and \cite{Yo}.   In the case of a $U(1)$-bundle over a compact surface, these were classified in \cite{St1}.   Ricci Yang-Mills solitons have also been studied in the context of dynamical systems.   In \cite{Jane}, it was discovered that the Ricci Yang-Mills flow is an ideal candidate for studying magnetic flows.  There is work in progress to determine whether Ricci Yang-Mills solitons have new dynamical properties.\\


This note is organized as follows.  In $\S 2$, we define both the Ricci Yang-Mills flow and the Ricci Yang-Mills soliton equations.  We provide the framework for the Ricci Yang-Mills soliton equations on 2-step nilpotent Lie groups in $\S 3$.  In this section we translate the notion of Ricci Yang-Mills solitons into the notion of distinguished points from Geometric Invariant Theory.  Finally in $\S 4$, we provide several examples (using algebraic techniques) of Lie groups that do not admit Ricci solitons but that do admit Ricci Yang-Mills solitons.\\

\textit{Acknowledgements.}  This note is a component of a larger project to understand the Ricci Yang-Mills flow and its special solutions.  The authors would like to thank Dan Jane for many enlightening conversations and Pat Eberlein for spotting a critical typo.

\section{Ricci Yang-Mills solitons}
%

The Ricci Yang-Mills flow  is a natural coupling of the Ricci flow and the Yang-Mills heat flow.   Let $\pi:P\to M$ be a principal bundle with structure group $G$.  Also, let $g$ be a metric on $M$, $k$ an $\mathfrak{Ad}$-invariant metric on $\mathfrak{g}$, and $\omega$ the connection 1-form on $P$.  We will consider so-called \emph{bundle metrics} on $P$ of the form
\[h=\pi^\ast g+k\omega,
\]
where $k \omega$ acts on vector fields by $k\omega(Y,Z)=k(\omega(Y),\omega(Z))$.

Writing the Ricci flow equations for a metric of this form with the additional hypothesis that the size of the fiber remains fixed, one can define the Ricci Yang-Mills flow to be
\begin{equation}
\label{rym total}
\frac{\partial h}{\partial t}=-2(Rc-Rc^V),
\end{equation}
where $Rc^V$ is the projection of the Ricci tensor onto its vertical component (cf.~\cite{St} and the proof of Theorem \ref{thm: RYM is invariant under bundle auts}).  If $G$ is abelian, using the definition of $h$ and the structure of $Rc(h)$, one can show that this equation is equivalent to the following system of equations:
\begin{subequations}
\label{rym}
\begin{align}
\frac{\partial g}{\partial t}&=-2Rc(g)+\tilde\Omega^2\\
\frac{\partial \tilde \omega}{\partial t}&=-\delta \tilde \Omega.
\end{align}
\end{subequations}
Here $\tilde \omega$ and $\tilde\Omega$ are the pullbacks under a local section of the connection 1-form and the bundle curvature, respectively.  Recall that when $G$ is abelian, $\tilde \Omega$ is a well-defined Lie algebra-valued 2-form on the base.  In coordinates, $\tilde \Omega^2_{ij}=g^{kl}k^{\alpha \beta}\tilde \Omega_{\alpha ki} \tilde{\Omega}_{\beta lj}$, where the greek indices are the Lie algebra indices and the Roman indices correspond to quantities measured with respect to $g$.   Existence and uniqueness of solutions to the Ricci Yang-Mills flow have been studied in \cite{St} and \cite{Yo}.\\

We would like to define Ricci Yang-Mills solitons in a way that is analogous to Ricci solitons.  Namely, due to the diffeomorphism invariance of the Ricci flow, one can define Ricci solitons to be fixed points of the Ricci flow in the space $\mathfrak{Met}/\mathfrak{Diff}$.  However, since a principal bundle is a manifold endowed with a group action of $G$, in defining Ricci Yang-Mills solitons, we are interested in diffeomorphisms that preserve the full structure of  the principal bundle.

\begin{defin}
An automorphism of a principal  bundle $\pi:P\to M$ is a diffeomorphism $F:P\to P$ such that $F(pg)=F(p)g$ for all $g\in G, p\in P$.  We denote this set by $\mathfrak{Diff}_G$.
\end{defin}

Let $F\in \mathfrak{Diff}_G$ be an automorphism of the principal bundle $P$.   Then $F$ descends to a diffeomorphism $f$ on $M$ satisfying  $f\circ \pi = \pi \circ F$, and  if $\omega$ is a connection on $P$, then $F^*\omega$ is also a connection.

\begin{thm}\label{thm: RYM is invariant under bundle auts}
The Ricci Yang-Mills flow is invariant under automorphisms of $P$.
\end{thm}

\begin{proof}

Using the  properties above, one immediately sees that $F^*h$ is a bundle metric for $F\in \mathfrak{Diff}_G$.  More precisely,
    $$F^*h = \pi^* (f^*g) + k F^*\omega.$$

For each $z\in \mathfrak{g}$, there is a canonical vector field on $P$ defined by $Z_p = \ddtat p\cdot exp(tz)$.  As $F$ preserves the $G$ action and $\omega (\ddtat p\cdot exp(tz) ) = z$, we see that $F_* (Z_p) = Z_{F(p)}$ and $\omega (Z_p) = \omega (Z_{F(p)})=(F^*\omega)(Z_p) = z$.  Observe that if $\{z_i\}$ is an orthonormal basis of $\mathfrak{g}$ then the induced vector fields $\{Z_i\}$ form an orthonormal frame of the vertical space relative to  both metrics, $h$ and $F^* h$.

We can write the Ricci Yang-Mills flow as $\frac{\partial h}{\partial t}=-2(Rc-Rc^V)$, where $Rc^V$ is the projection of the Ricci tensor onto its vertical component.   Specifically, if $U$ is a vector field on $P$, then we can define the projection onto its vertical component to be $p_h(U)=\sum_i h(Z_i,U) Z_i$, where $\{Z_i\}$ is an orthonormal basis (relative to $h$) of the vertical space as above.  In this notation, we have  $Rc^V(h)(U,V)=Rc(p_h(U),p_h(V))$.

Let $F:P\to P$ be a bundle automorphism.  Since $F$ is a diffeomorphism, clearly $Rc(F^\ast h)=F^\ast Rc(h)$.  It remains only to check that $Rc^V(F^\ast h)=F^\ast Rc^V(h)$.  First we show that $p_{F^\ast h} = p_h \circ F_* $.  By definition,
\begin{eqnarray*}
p_{F^\ast h}(U)&=& \sum_i (F^* h)(Z_i,U) Z_i\\
		       &=& \sum_i (f^* g) (\pi_* Z_i, \pi_* U) + k(
                \omega(F_*U),\omega(F_*Z_i))Z_i\\
		       &=& \sum_i k( \omega(F_*U),\omega(Z_i) )Z_i\\
                &=& \sum_i h(F_* U, Z_i) Z_i\\
                &=& p_h ( F_* U)
\end{eqnarray*}
Here we have used the fact that $\{Z_i\}$ will be orthonormal in both metrics $h$ and $F^*h$.  Thus $p_{F^\ast h} = p_h \circ F_\ast $.  Using this fact and the diffeomorphism invariance of $Rc$, one sees that in fact $Rc^V(F^\ast h)=F^\ast Rc^V(h)$.  Thus the Ricci Yang-Mills flow is invariant under bundle automorphisms of $P$.
\end{proof}

\begin{cor} The Ricci Yang-Mills flow preserves the set of left-invariant metrics on a Lie group $N$.
\end{cor}
\begin{proof}
Left multiplication $L_g(p)=gp$ is a bundle automorphism since left and right multiplication commute.  Thus the result follows from above.
\end{proof}

Since the set of left-invariant metrics is preserved under this evolution, we can interpret  the Ricci Yang-Mills flow as an evolution of the metric on a single tangent space; more precisely, we will evolve an inner product on the Lie algebra.  This is the standard approach to studying geometric evolutions on Lie groups with left-invariant metrics.

\subsection{Self-similar solutions to Ricci Yang-Mills equations.}

Analogous to the case of Ricci flow, we define Ricci Yang-Mills solitons to be generalized fixed points of Eq.~\ref{rym}.

\begin{defin}
A solution  $(g_t, \tilde\omega_t)$ to the Ricci Yang-Mills equations is a self similar solution if there exists a scaling $\sigma(t)$ and a family of diffeomorphisms $\psi_t$ $ \in \mathfrak{Diff}(M)$ such that $g_t = \sigma(t) \cdot \psi_t^* g$ and $\tilde\Omega_t = \psi_t^* \tilde\Omega$ (with $\sigma(0) = 1$ and $\psi_0$ the identity). Let $X \in \Gamma(M, TM)$ generate $\psi$ near $t = 0$.  As in the case of Ricci flow, one can show that the notion of self-similar solutions is equivalent to $(g,\tilde\omega)$ satisfying
\begin{subequations}
\label{RYM soliton}
\begin{align}
(\tilde\Omega^2-2Rc)(g) &= \sigma'(0) g + \mathcal{L}_X g,\\
 \Delta_d \tilde\Omega &= \mathcal{L}_X \tilde\Omega,
\end{align}
\end{subequations}
where $\Delta_d$ is the Hodge Laplacian.  We will call solutions satisfying Eq.~\ref{RYM soliton}  \emph{(weak) Ricci Yang-Mills solitons}.
\end{defin}

\textit{Notation:} Let $\lambda = \sigma '(0)$.  If $\lambda>0$ we say the Ricci Yang-Mills soliton is an expander, if $\lambda=0$ the Ricci Yang-Mills soliton is called  steady, and $\lambda<0$ is called a shrinker.

\begin{defin}\label{defin: strong RYM soliton}  In contrast to weak Ricci Yang-Mills solitons, we say that a Ricci Yang-Mills soliton is a strong Ricci Yang-Mills soliton when there exists a family $\tilde \psi_t \in \mathfrak{Diff}_G$ which induces $\psi_t$; that is, such that $\pi \circ \psi_t = \varphi_t \circ \pi$.
\end{defin}

\noindent
\begin{remark}
The notion of being a weak Ricci Yang-Mills soliton is intrinsic to the base manifold; as the group $G$ is abelian, we may consider $g$ and $\tilde\Omega$  as objects living on $M$.  From the perspective of the base, one doesn't see diffeomorphisms of the total space.   It is not known, even in the case of $U(1)$-bundles over surfaces, whether weak solitons are always strong.  Considering not all diffeomorphisms on $M$ are induced by bundle automorphisms of $P\to M$, we expect there to exist Ricci Yang-Mills solitons which are precisely weak.  This will be the object of future study.
\end{remark}

In the setting of nilpotent Lie groups, our principal bundles are $N \to N/Z$ where $Z$ is the center.  Here every diffeomorphism lifts to a diffeomorphism of the total space and hence all Ricci Yang-Mills solitons will be strong in this paper, see Section \ref{section: RYM sol of Lie type}.\\

\noindent
\textit{Remark.}
An Einstein Yang-Mills metric as defined in \cite{St} and \cite{Yo} is one such that the metric on the base is Einstein and the connection is Yang-Mills; i.e. $\delta \tilde\Omega=0$.  On a 2-step nilpotent Lie group, every metric is a metric of this type (cf. Lemmas \ref{lemma: base is flat} and \ref{lemma: Delta Omega = 0}).

Ricci Yang-Mills solitons are not direct generalizations of Einstein Yang-Mills metrics (as defined above) in the same way that Ricci solitons are generalizations of Einstein metrics.    Recall that Ricci solitons are fixed points of the volume-normalized Ricci flow, which differs from the Ricci flow only by a change of scale in space and time.  Einstein Yang-Mills metrics (as defined above) are fixed points of a certain volume normalized Ricci Yang-Mills flow; however, due to the lack of  scale invariance of this equation, the volume normalized flow does not differ only by a change of scale (see \cite{St} or \cite{Yo}).

An alternate definition of Einstein Yang-Mills, which is more natural from the view point of special solutions to the Ricci Yang-Mills flow, would be a Ricci Yang-Mills soliton generated by a trivial vector field; that is, $L_X=0$.  This is consistent with the idea of an Einstein metric from the perspective of Ricci solitons and Ricci flow.  To avoid confusion with previous definitions, we make the following definition
\begin{defin} A \emph{trivial Ricci Yang-Mills soliton} is one that is generated  by the trivial vector field $X=0$.\end{defin}

 Below we will construct many examples of trivial Ricci Yang-Mills solitons; these correspond to so-called minimal points of a particular representation (see Corollary \ref{cor: closed SL q C orbit implies RYM soliton}).

\section{Ricci Yang-Mills solitons on nilmanifolds}

A natural test case in the search for Ricci Yang-Mills solitons on principal bundles is the case of a  torus bundle over a  torus.  These compact manifolds are precisely the locally homogeneous manifolds which are modeled on 2-step nilpotent Lie groups (see \cite{Pal61}).  More precisely, these spaces are quotients of nilpotent Lie groups by cocompact lattices.  As in the case of Ricci flow, to understand the dynamics of this geometric evolution on a compact manifold, we study the evolution on the simply connected cover, a nilpotent Lie group with a left-invariant metric.  This cover is also a principal bundle and the covering map is a morphism of bundles.  We recall some basic facts for the convenience of the reader.

\begin{defin}Let $\mathfrak{N}$ be a finite dimensional Lie algebra, and for $i\geq 1$, let $\mathfrak{N}^i=[\mathfrak{N},\mathfrak{N}^{i-1}]$, where $\mathfrak{N}^0=\mathfrak{N}$.  Then $\mathfrak{N}$ is said to be nilpotent if $\mathfrak{N}^i=\{0\}$ for some $i$.  A nilpotent Lie algebra is $k$-step if $\mathfrak{N}^k=\{0\}$ but $\mathfrak{N}^{k-1}\neq \{0\}$.  A Lie group is said to be ($k$-step)  nilpotent if its Lie algebra is ($k$-step) nilpotent.
\end{defin}

\begin{defin}\label{defin: type p,q} A 2-step nilpotent Lie group $N$ or Lie algebra $\mathfrak N$ is said to be of type $(p,q)$ if dim $[\mathfrak N,\mathfrak N] = p$ and codim $[\mathfrak N,\mathfrak N] =q$.
\end{defin}

Observe that $p$ above satisfies $1\leq p \leq \frac{1}{2}q(q-1) = \dim \mathfrak{so}(q)$.  Stratifying the space of 2-step nilpotent Lie algebras into types $(p,q)$ is very convenient in terms of phrasing generic results.

\begin{defin} Let $\{X_i\}$ be a basis of the Lie algebra $\mathfrak N$.  The structure constants relative to the basis $\{X_i\}$ are the coefficients $\{c_{ij}^k\}$ defined by $[X_i,X_j]=\sum_k c_{ij}^k X_k$.  \end{defin}

\begin{defin}  Let $N$ be a Lie group with Lie algebra $\mathfrak N$.  There exists a map $exp_N : \mathfrak N \to N$  called the Lie group exponential.  When there is no confusion, we write $exp$ for $exp_N$.
\end{defin}

Recall that when $N$ is a simply connected nilpotent Lie group, the map $exp_N$ is a diffeomorphism from $\mathfrak N$ to $N$.

\begin{thm}[Campbell-Baker-Hausdorff formula]\label{thm: CBH formula} For $X,Y \in \mathfrak N$, we have
    $$exp(X) \ exp(Y) = exp( X+Y + \frac{1}{2} [X,Y] + \ higher \ order \ terms )$$
where `higher order terms' means combinations of 2 or more brackets involving $X$ and $Y$.
\end{thm}

This formula explicitly relates the Lie product of the group and the Lie bracket of the algebra.  The higher order terms can be written explicitly, but we will only use this formula in the case that one of $X$ or $Y$ is in the center; in this case we have $exp(X)exp(Y)=exp(X+Y)$ (see Section \ref{section: RYM sol of Lie type}).\\

Every nilpotent Lie group $N$ can be viewed as a non-trivial principal bundle.  The total space will be $P=N$, and the group $G$ will be the center of $N$ acting on the right.  We describe this in detail in the 2-step case.  The details in the general case are similar.

\subsection{2-step nilmanifolds}

We would like to consider the Ricci Yang-Mills soliton equations on simply-connected 2-step nilpotent Lie groups.  Let $N$ be a 2-step nilpotent Lie group with Lie algebra $\mathfrak N$.  We endow $N$ with a left-invariant metric $h$; this is equivalent to endowing the Lie algebra $\mathfrak N$ with an inner product.  Let $Z=[N,N]$ be the commutator subgroup, $\mathfrak Z = Lie \ Z$, and $\mathcal H = \mathfrak Z ^ \perp$ the orthogonal complement relative to the given metric on $\mathfrak N$.  We point out that $Z$ is a central subgroup since $N$ is 2-step nilpotent.  One could choose to work with either the full center of $N$ or the commutator $[N,N]$ in what follows.

Let $\{X_1,\dots, X_q\} \cup \{Z_1,\dots, Z_p\}$ be an orthonormal basis of $\mathfrak N = \mathcal H \oplus \mathfrak Z$.  Here $q=\dim \mathcal H$, $p=\dim \mathfrak Z$, and $n = q+p = \dim N$; this 2-step nilpotent algebra is of type $(p,q)$. By left-translating, we can treat this basis of $\mathfrak N$ as a left-invariant frame on $N$.  Relative to this basis we have the (Lie algebra) structure coefficients defined via
    $$ [ X_i, X_j] = \sum_k c_{ij}^k Z_k$$.

Notice that $\pi:N\to N/Z$ is naturally a principal $G$-bundle where $G = \mathbb R^p \simeq Z$.  The action of $G\simeq Z$ will be given by first injecting $Z$ into $N$ and then multiplying on the right.  To distinguish between $\mathfrak Z$ abstractly versus embedded in $\mathfrak N$, we will use lower case letters to denote elements of $\mathfrak g$ and upper case letters to denote elements of $\mathfrak Z\subset \mathfrak N$; that is, given $z\in \mathfrak g$, $z\to Z \in \mathfrak Z \subset \mathfrak N$.  A  choice of a horizontal subspace $\mathcal{H}$ yields a connection 1-form $\omega$ which vanishes on $\mathcal H$  and takes values in the Lie algebra $\mathfrak{g}$.  More precisely, we define our connection so that $\omega(Z)=z$ and $\omega(X)=0$ for $Z\in \mathfrak Z$, $X \in \mathcal H$.  Thus far, we have defined our connection on $\mathfrak N = T_e N$. We extend the definition of the connection to the rest of $N$ by imposing $\omega$ be left-invariant under $N$.

Recall that the Ricci Yang-Mills soliton equations are
\begin{subequations}
\label{rym 1}
\begin{align}
-2Rc_g+\tilde\Omega^2&=\mathcal{L}_Xg+\lambda g\\
\Delta_d\tilde\Omega&=\mathcal{L}_X\tilde\Omega.
\end{align}
\end{subequations}
Here $g$ is the induced metric on $M=N/Z$.  As $Z$ is a normal subgroup of $N$, $N/Z$ is a Lie group and the metric $g$ is left $N/Z$-invariant.\\

\begin{lemma}\label{lemma: base is flat} In the 2-step nilpotent setting, Eq.~\ref{rym 1}a becomes $\tilde\Omega^2=\mathcal{L}_Xg+\lambda g$.
\end{lemma}
\begin{proof}
 As the base is an abelian Lie group $N/Z$ with left-invariant metric, it is flat and hence $Rc_g=0$.
\end{proof}

\begin{lemma}\label{lemma: Omega vs structure matrices} Let $X_i, X_j$ be horizontal vectors in the basis of $\mathfrak N$ above, then
$\Omega (X_i,X_j) = \sum_\alpha -c_{ij}^\alpha z_\alpha$.
\end{lemma}
\begin{proof}
By definition of $\Omega$,
\begin{displaymath}
\Omega(X_i, X_j) = d \omega (X_i, X_j) = X_i( \omega (X_j)) - X_j( \omega(X_i)) - \omega ([X_i, X_j])=-\omega(\sum_k c_{ij}^kZ_k)=\sum_k -c_{ij}^kz_k.
\end{displaymath}

Here we have used that fact that $\omega$ is left-invariant and $X_i$ is a left-invariant vector field on $N$.  Thus, $\omega(X_i)$ is a constant function on $N$ and $X_j(\omega(X_i))=0$.
\end{proof}

\begin{lemma}\label{lemma: Delta Omega = 0} The connection $\tilde \omega$ is Yang-Mills; that is,
$\Delta_d \tilde{\Omega}=0$.
\end{lemma}
\begin{proof}
Recall that $\tilde\Omega = d\tilde\omega$ and so $\delta \tilde\Omega = 0$ if and only if $\Delta_d \tilde\Omega = 0$.  Let $U \subset M$ be an open set, and let $s:U\to N$ be a local section.  We can define $\tilde{X}_i=\pi_\ast X_i$ to be a left-$N/Z$-invariant vector field on the base.   We have that $[X_i,X_j]=c_{ij}^kZ_k$.  Thus
\begin{eqnarray*}
\tilde{\Omega}(\tilde{X}_i,\tilde{X}_j)&=&s^\ast \Omega(\tilde{X}_i,\tilde{X}_j)\\
							&=&\Omega(s_\ast\tilde{X}_i,s_\ast\tilde{X}_j)\\
							&=&\Omega(X_i,X_j).
\end{eqnarray*}
So $\tilde{\Omega}=\sum_{ijk}-c_{ij}^k z_k\tilde \sigma^i \wedge \tilde \sigma^j$, where $\{\tilde \sigma^i\}$ are dual to $\{\tilde X_i \}$.  To compute $\Delta_d\tilde{\Omega}$, we only need to compute $d\star d\star (\tilde{\sigma}^i\wedge \tilde{\sigma}^j)$, as the structure constants come out.  We compute that $d\star(\tilde{\sigma}^i\wedge \tilde{\sigma}^j)=(-1)^{i+j-1}d(\tilde{\sigma}^1\wedge\cdots \wedge \hat{i}\hat{j} \wedge \cdots \wedge \tilde{\sigma}^q)$, and
\[
d\tilde{\sigma}^k(\tilde{X}_i,\tilde{X}_j)=\tilde{X}_i(\tilde{\sigma}^k(\tilde{X}_j)-\tilde{X}_j(\tilde{\sigma}^k(\tilde{X}_i)-\tilde{\sigma}^k([\tilde{X}_i,\tilde{X}_j])=-\tilde{\sigma}^k([\tilde{X}_i,\tilde{X}_j]).
\]
Again, we have used left-invariance to make two of the middle terms vanish in the above equation.  Lastly,  $[\tilde{X}_i,\tilde{X}_j]=0$ as $\mathfrak N / \mathfrak Z$ is abelian.

Thus $\Delta_d\tilde{\Omega}=0$.
\end{proof}

\begin{lemma} In the 2-step nilpotent setting, Eq.~\ref{rym 1}b becomes $\mathcal{L}_X\tilde{\Omega}=0$.
\end{lemma}

\begin{prop}\label{prop: RYM from deriv on N/Z}
Let $D$ be a derivation of $\mathfrak N/ \mathfrak Z$ and $exp(tD)$ the associated 1-parameter family of automorphisms of $N/Z$.  Using this family of diffeomorphisms, the Ricci Yang-Mills soliton equations on 2-step nilpotent Lie groups become
\begin{subequations}
\label{rym w/deriv}
\begin{align}
\tilde \Omega^2(v,w)&=g((D+D^t)v,w)+\lambda g(v,w)\\
0&=\tilde \Omega(Dv,w)+\tilde \Omega(v,Dw),
\end{align}
\end{subequations}
where $v, w$ are vector fields on $M = N/Z$.
\end{prop}
\begin{proof}  We postpone presenting the details of this proof.  See Appendix A for information regarding one-parameter families of automorphisms.
\end{proof}

\begin{observ} Every linear map of an abelian Lie algebra is an automorphism of the Lie algebra.
\end{observ}

\noindent Thus, $\mathfrak N$ being a 2-step nilpotent Lie algebra means $\mathfrak N/\mathfrak Z$ is abelian  and we can use any linear map ${ D: \mathfrak N / \mathfrak Z \to \mathfrak N / \mathfrak Z}$ in Eq.~\ref{rym w/deriv}.

\begin{defin}\label{defin: RYM soliton of Lie type} We say that a left-invariant Ricci Yang-Mills soliton is of \emph{Lie type} if it comes from an automorphism of the base $N/Z$ as in Proposition \ref{prop: RYM from deriv on N/Z}.  We say that such a metric is of \emph{symmetric Lie type} if the derivation is symmetric with respect to the inner product on $\mathfrak N/ \mathfrak Z$.
\end{defin}

In the sequel, we build many examples of such Ricci Yang-Mills solitons. Presently, we have no examples of left-invariant Ricci Yang-Mills solitons which are not of symmetric Lie type.  This question will be investigated in future work.

\subsection{2-step nilpotent from the perspective of structure matrices}
In this section we approach this problem of finding Ricci Yang-Mills solitons using so-called `structure matrices'.  Studying 2-step nilpotent Lie groups from this point of view is quite natural and has been used by Eberlein, Jablonski, and others to obtain results about Ricci solitons.  Using this approach we prove the existence of a large class of Ricci Yang-Mills solitons using Geometric Invariant Theory (cf. Theorem \ref{thm: RYM soliton vs distinguished pt}) and construct examples of nilmanifolds that do not admit left-invariant Ricci solitons but do admit left-invariant Ricci Yang-Mills solitons (see Section \ref{section: examples}).\\

The geometry of $N$ can be completely encoded  by a tuple of structure matrices $(C^1,\dots, C^p)$ which is an element of $\mathfrak{so}(q)^p$.  We establish this perspective below.  Our main references are \cite{Eber07} and \cite{Jab:Thesis}.

Recall that a nilpotent Lie group $N$ with left-invariant metric $\langle,\rangle $ is equivalent to a nilpotent Lie algebra $\mathfrak N$ with inner product, also denoted $\langle,\rangle $.  Let $\{X_1,\dots,X_q\} \cup \{Z_1,\dots ,Z_p\}$ be an orthonormal basis of $\mathfrak N$ where the $Z_k$ form an orthonormal basis of $\mathfrak Z$.  Relative to this basis we may compute the structure constants $\{c_{ij}^k\}$ defined via $[X_i,X_j] = \sum_k c_{ij}^k Z_k$.  Thus we may associate to our basis a $p$-tuple of matrices $(C^1,\dots,C^p)$ where $(C^k)_{ij}=c_{ij}^k$.  Notice that different $p$-tuples of matrices can be associated to a given  $\mathfrak N$.   We  describe below how these different tuples of matrices are related to each other.

Conversely, given a tuple $C=(C^1,\dots,C^p)$, we can naturally associate to it a metric 2-step nilpotent Lie algebra.  This construction is dual to the construction of $p$-tuples above.   We will require the $C^k$ to be linearly independent in $\mathfrak{so}(q)$ so that the commutator of the constructed nilpotent algebra will have dimension $p$.

Let $\{e_1,\dots,e_q,e_{q+1},\dots,e_{q+p} \}$ be the standard basis of $\mathbb R^q\oplus \mathbb R^p$.  Endow $\mathbb R^q\oplus \mathbb R^p$ with the standard inner product so that this basis is orthonormal.  Define the Lie bracket as
    $$[e_i,e_j] = \sum_k C^k_{ij}e_{q+k}$$
for $1\leq i,j \leq q$, and let all other brackets be trivial.  This Lie algebra is clearly 2-step nilpotent with commutator equal to $\mathbb R^p = span\langle e_{q+1}, \dots , e_{q+p}\rangle $.  We denote the metric 2-step nilpotent Lie algebra associated to $C$ by $\mathfrak N_C$ and the 2-step nilpotent Lie group with left-invariant metric associated to $\mathfrak N_C$ by $N_C$.

Using tuples of matrices, we may study 2-step nilpotent Lie algebras not just individually, but as elements of the much larger space  $\mathfrak{so}(q)^p$.  As stated above, elements of $\mathfrak{so}(q)^p$ whose coordinates are linearly independent correspond to 2-step (metric) nilpotent  Lie algebras of type $(p,q)$ (cf. Definition \ref{defin: type p,q}).

\subsection{Geometric Invariant Theory and 2-step nilpotent Lie algebras}
On the vector space $\mathfrak{so}(q)^p$ there is a natural action of $GL_q\mathbb R \times GL_p\mathbb R$ which is a linear representation.  Given $g\in GL_q\mathbb R$ and $C=(C^1,\dots, C^p)$, we define
    $$ g\cdot C = (gC^1g^t, \dots, gC^pg^t)$$
which acts on each coordinate individually.  It is clear that $g\cdot C^k = gC^kg^t \in \mathfrak{so}(q)$ given that $C^k\in \mathfrak{so}(q)$.  The action of $GL_p\mathbb R$ takes linear combinations of the coordinates of $C=(C^1,\dots ,C^p)$.  Given $h\in GL_p\mathbb R$ we define
    $$h\cdot C = D=(D^1,\dots, D^p) \quad \mbox{ with } \quad D^k = \sum_l h_{lk}C^l$$
One can compute that these actions commute and thus we have an action of $GL_q\mathbb R \times GL_p\mathbb R$ on $\mathfrak{so}(q)^p$.

Using this action, we can easily describe  when two different structure matrices produce the same nilpotent Lie group.  For proofs of the next two theorems we refer the reader to \cite{Eber07}.

\begin{thm} Let $C, D \in \mathfrak{so}(q)^p$ correspond to 2-step nilpotent Lie algebras $N_C, N_D$, respectively.  Then $N_C$ and $N_D$ are isomorphic Lie groups if and only if $D \in GL_q\mathbb R \times GL_p \mathbb R \cdot C$, the orbit of $C$ under the group action of $GL_q\mathbb R \times GL_p\mathbb R$.
\end{thm}
Here we were only concerned with the underlying Lie group structure of $N_C$ and $N_D$.  The next theorem considers the metric structures as well.

\begin{thm}  Let $C, D \in \mathfrak{so}(q)^p$ correspond to 2-step nilpotent Lie algebras $N_C, N_D$, respectively, with left-invariant metrics.  Then $N_C$ and $N_D$ are isometric as Riemannian manifolds if and only if $D \in O(q)\times O(p) \cdot C$, the orbit of $C$ under the compact subgroup $O(q)\times O(p) \subset GL_q\mathbb R \times GL_p\mathbb R$.
\end{thm}

This representation of $GL_q\mathbb R\times GL_p\mathbb R $ on $\mathfrak{so}(q)^p$ has even more structure from the view point of Geometric Invariant Theory.  Once translated, these extra structures on the representation space have strong consequences on the Riemannian geometry of associated nilpotent Lie groups.  We present a brief discussion below and refer the reader to \cite{Jab:DetectingOrbits} for a more thorough treatment.\\

Associated to the representation of $GL_q\mathbb R\times GL_p\mathbb R$ on $\mathfrak{so}(q)^p$, we have a Lie algebra representation of $\mathfrak{gl}_q\mathbb R \times \mathfrak{gl}_p\mathbb R$ on $\mathfrak{so}(q)^p$.  This is obtained in the usual way via differentiation; that is, given $(X,Y)\in \mathfrak{gl}_q\mathbb R \times \mathfrak{gl}_p\mathbb R$ and $C\in \mathfrak{so}(q)^p$ we have
    $$ (X,Y)\cdot C = X\cdot C + Y\cdot C,$$
where $X\cdot C = (XC^1+C^1X^t, \dots , XC^p + C^pX^t)$ and $Y\cdot C = D$ with $D^k = \sum_l Y_{lk} C^l$.

The space $\mathfrak{so}(q)$ has the inner product $\langle C,D\rangle  = tr(CD^t) = - tr (CD)$.  This inner product extends to $\mathfrak{so}(q)^p$ by making the coordinates of the tuple orthogonal; that is, consider $C=(C^1,\dots,C^p)$ and $D=(D^1,\dots,D^p)$ then
    $$\langle C,D\rangle  =  \langle (C^1,\dots,C^p),(D^1,\dots,D^p)\rangle  = \sum_\alpha \langle C^\alpha,D^\alpha\rangle  = \sum_\alpha -tr(C^\alpha D^\alpha).$$
We define two polynomials, $m_1$ and $m_2$, which are associated to our representation of $GL_q\mathbb R \times GL_p\mathbb R$.   Here $m_1$ is the moment map for the action of $GL(q,\mathbb{R})$ on $\mathfrak{so}(q)^p$, and $m_2$ is the moment map for the action of $GL(p,\mathbb{R})$ on $\mathfrak{so}(q)^p$.  Notice that the polynomial $m_1$ will be valued in $symm_q$, the symmetric $q\times q$ matrices, while the polynomial $m_2$ will be valued in $symm_p$.  For $C\in \mathfrak{so}(q)^p$, we define
    \begin{eqnarray*} m_1(C) &=& -2 \sum_\alpha (C^\alpha)^2 \\
        m_2(C)_{ij} &=&  \langle C^i,C^j\rangle .
    \end{eqnarray*}
Adding these together one has the moment map $m=m_1+m_2$ for the action of $GL_q\mathbb R \times GL_p\mathbb R$.  This is valued in $symm_q \oplus symm_p$.

We are interested in three different group actions on $\mathfrak{so}(q)^p$; namely, the actions of the full group $GL_q\mathbb R\times GL_p\mathbb R$ and its  subgroups  $GL_q\mathbb R$, the first factor, and $SL_q\mathbb R \subset GL_q\mathbb R$.  In the following definition, $G$ will denote one of these three groups, and $m_G$ will denote the moment map corresponding to $G$.

\begin{defin}\label{defin: distinguished points} We call a point $C\in \mathfrak{so}(q)^p$ $G$-distinguished if $m_G(C)\cdot C = r C$ for some $r\in \mathbb R$.  We call a point $G$-minimal if $m_G(C)\cdot C =0$.  Minimal points are obviously a special kind of distinguished point.
\end{defin}

Here $m_G(C) \in \mathfrak{gl}_q\mathbb R \times \mathfrak{gl}_p\mathbb R$, since $G$ is a subgroup of $GL_q\mathbb R \times GL_p\mathbb R$, and it acts via the Lie algebra action of $\mathfrak{gl}_q\mathbb R \times \mathfrak{gl}_p\mathbb R$.  Distinguished and minimal points can be defined more generally for any representation of a reductive group on a vector space \cite{Jab:DistinguishedOrbits}.

\begin{prop} Let $C$ be a distinguished point as above, then $r\geq 0$.
\end{prop}

This is a consequence of the facts that $\nabla ||m_G||^2  (C) = m_G(C)\cdot C$ and that a function is non-decreasing along its gradient flow; for more details see \cite{Jab:DistinguishedOrbits}.  In this setting distinguished points play a very interesting role (cf. Theorems \ref{thm: Ricci solion vs. distinguished point}, \ref{thm: geod flow inv vs. distinguished pt}, and \ref{thm: RYM soliton vs distinguished pt}).

\begin{thm}\label{thm: Ricci solion vs. distinguished point} Let $C\in \mathfrak{so}(q)^p$ correspond to a 2-step nilpotent Lie group $N_C$ with left-invariant metric.  Then $N_C$ is a nilsoliton if and only if $C$ is a distinguished point of the $GL_q\mathbb R\times GL_p\mathbb R$ action.
\end{thm}
By nilsoliton we mean a nilpotent Lie group with left-invariant Ricci soliton metric.  This was originally proven for all nilpotent Lie groups (not just 2-step) by Jorge Lauret \cite{La}.  In the 2-step nilpotent setting, Eberlein \cite{Eber07} proved this in the frame work of structure matrices; this approach has proven very fruitful for constructing examples.

\begin{thm}\label{thm: geod flow inv vs. distinguished pt} Let $C\in \mathfrak{so}(q)^p$ corresponding to a 2-step nilpotent Lie group $N_C$ with left-invariant metric.  Then the metric on $N_C$ is so-called geodesically flow invariant if and only if $C$ is a minimal point of the $SL_q \mathbb R$ action.
\end{thm}

We have not defined the notion of a metric being `geodesically flow invariant' and refer the reader to \cite{Eber07} for more details.  This theorem is presented so that the reader may place Theorem \ref{thm: RYM soliton vs distinguished pt} in a more general context with Theorems \ref{thm: Ricci solion vs. distinguished point} and \ref{thm: geod flow inv vs. distinguished pt}.

\begin{question}\label{question: distinguished pts of GLqR action} Is there good geometric meaning to Lie groups $N_C$ which correspond to a point $C$ which is a distinguished point of the $GL_q\mathbb R$-action?
\end{question}

This question was asked to us by Pat Eberlein and has been answered in Theorem \ref{thm: RYM soliton vs distinguished pt}.

\subsection{Ricci Yang-Mills equations and Geometric Invariant Theory}
As Ricci solitons are naturally associated to distinguished points (cf. Definition \ref{defin: distinguished points} and Theorem \ref{thm: Ricci solion vs. distinguished point}), one is naturally lead to investigate if there is a similar connection between Ricci Yang-Mills solitons and Geometric Invariant Theory.  There is a strong (and similar) relationship in the Ricci Yang-Mills setting.  We will study the Ricci Yang-Mills soliton equations from the perspective of structure matrices.  We begin by translating   Eq.~\ref{rym w/deriv} into a system of equations on tuples of matrices.

In the following proposition, $C\in \mathfrak{so}(q)^p$ corresponds to the metric Lie group $N_C$, and $D$ is an element of $\mathfrak{gl}_q\mathbb R$.  Additionally, we present the Ricci Yang-Mills equations here as $(1,1)$ tensors as opposed to $(2,0)$ tensors.  From the perspective of structure matrices, it is more natural to present the equations in the following form.
\begin{prop}
The Ricci Yang-Mills soliton equations of Lie type (Eqns.  \ref{rym w/deriv} a \& b) on 2-step nilpotent Lie groups can be written as
    \begin{subequations}
    \label{rym 2}
    \begin{align}
    m_1(C) = -2\sum_\alpha (C^\alpha)^2  &= 2\lambda Id + 2 (D+D^t)\\
            0 &= D^t \cdot C
    \end{align}
    \end{subequations}
where $D^t \cdot C$ denotes the action of $\mathfrak{gl}_q$ on $\mathfrak{so}(q)^p$; that is, $D^t \cdot C = (D^t \cdot C^1, \dots , D^t \cdot C^p)$ and $D^t \cdot C^i = D^t C^i+C^iD$.
\end{prop}

\begin{proof}

We begin by raising an index on $\tilde \Omega^2$ and show that this is $\frac{1}{2}m_1(C)$.  Recall that \[\tilde \Omega^2_{ij} = \sum_{\alpha \beta k l} g^{kl}k^{\alpha \beta} \tilde \Omega_{\alpha k i} \tilde \Omega_{\beta l j}.\]
We will use the orthonormal basis $\{X_i\} \cup \{Z_\alpha\}$ of $\mathfrak N$ that was used to calculate our structure matrix.  Thus we have
    \begin{eqnarray*}
        \tilde \Omega^2_{ij} = \sum_{\alpha \beta k l} g^{kl}k^{\alpha \beta} \tilde \Omega_{\alpha k i} \tilde \Omega_{\beta l j} &=& \sum_{\alpha l} \tilde \Omega_{\alpha l i} \tilde \Omega_{\alpha l j}\\
            &=& \sum_{\alpha l} (-C^\alpha)_{ l i} (-C^\alpha)_{ l j} = - \sum_{\alpha} (C^\alpha)^2_{ij} = \frac{1}{2} m_1(C)_{ij}
    \end{eqnarray*}
Here we have used Lemma \ref{lemma: Omega vs structure matrices} to compare $\tilde \Omega ^2$ and $C$.  The right-hand side of Eq.~\ref{rym w/deriv}a is easily converted to a $(1,1)$-tensor to obtain the claimed result.

For the second equation, recall that $\tilde \Omega = - \sum z_k c_{ij}^k \tilde \sigma_i \wedge \tilde \sigma_j$ and $c_{ij}^k = \langle C^k X_i,X_j\rangle $.  Thus
\[\tilde \Omega(v,w) = - \sum_k z_k \langle C^kv,w\rangle ,\]
where we are identifying $\mathcal H \simeq \mathfrak N/ \mathfrak Z$ isometrically via $\pi_*$.  Therefore
\[\tilde \Omega(Dv,w) + \tilde \Omega(v,Dw) = - \sum_k z_k \langle C^kDv,w\rangle  + \langle C^kv,Dw\rangle  = - \sum_k z_k \langle (C^k D + D^t C^k) v,w \rangle \]
as required.
\end{proof}

\begin{thm}
\label{thm: RYM soliton vs distinguished pt}
Let $N_C$ be the metric 2-step nilpotent Lie group corresponding to a tuple $C\in \mathfrak{so}(q)^p$.  Then the metric on $N$ is a left-invariant Ricci Yang-Mills soliton of symmetric Lie type (cf. Definition \ref{defin: RYM soliton of Lie type}) if and only if $C$ is a distinguished point of the action of $GL_q\mathbb R$ on $\mathfrak{so}(q)^p$ (cf. Definition \ref{defin: distinguished points}).
\end{thm}

This theorem should be compared to its analogue for Ricci solitons (Theorem  \ref{thm: Ricci solion vs. distinguished point}).  An interesting and easily proved consequence of the above theorem is the following.

\begin{cor} Ricci Yang-Mills solitons of symmetric Lie type are expanders; i.e. have $\lambda >0$.
\end{cor}

\begin{proof}[Proof of theorem and corollary]  By definition, $C$ being a distinguished point of the action of $GL_q$ on $\mathfrak{so}(q)^p$ is equivalent to $m_1(C)\cdot C = a(C) C$ where $a(C) > 0$.  This holds if and only if $m_1(C)= \frac{1}{2} a(C) + B$ where $B\in Stab_C$, that is, $B\cdot C = 0$.  Since $m_1(C)$ is always a symmetric matrix, $B$ is a symmetric matrix.

Using Eq.~\ref{rym 2}, we see that if $D$ is symmetric, we have our equivalence using $a=2\lambda$ and $D=-4B$.  Lastly,  $a>0$ implies $\lambda >0$.
\end{proof}

\begin{remark}
Notice that $\lambda >0$ agrees with the sign convention of \cite{La08}.  There for nilpotent Lie groups, he defines the Ricci soliton equation to be $Rc=cI+D$ and shows that $c<0$.  Our first soliton equation has that $Rc^H=-\lambda I+D$, where $Rc^H$ is the horizontal component of the Ricci tensor. \end{remark}

At this point we are able to use general theorems from Geometric Invariant Theory to prove results about the existence of Ricci Yang-Mills solitons.

\begin{cor}\label{cor: closed SL q C orbit implies RYM soliton}  If the orbit $SL_q\mathbb R \cdot C$ is closed in $\mathfrak{so}(q)^p$ then the corresponding Lie group admits a (trivial) Ricci Yang-Mills soliton.  Moreover, such metrics are precisely the geodesically flow-invariant metrics (cf. Theorem \ref{thm: geod flow inv vs. distinguished pt}).
\end{cor}

This corollary follows from the fact that if $SL_q\mathbb R \cdot C$ is closed, then there exists a minimal point on the orbit (assume it is $C$) satisfying $m_1(C) = r\ Id$. See \cite{Jab:Thesis} for more details.  Once $m_1(C)$ has this form, it is clear that $C$ will be a (trivial) Ricci Yang-Mills soliton.  This corollary provides us with a very general procedure for building examples of (trivial) Ricci Yang-Mills solitons.

\begin{remark}  Consider 2-step nilpotent Lie groups  of type $(p,q)$.  If $ p \leq \frac{1}{2}q(q-1) -2$ then almost every  $SL_q\mathbb R$-orbit is  closed \cite{Jab:Thesis}.  Hence, almost every 2-step nilpotent Lie group admits a Ricci Yang-Mills soliton when $p \leq \frac{1}{2}q(q-1) -2$.  \end{remark}

In the nilpotent setting, we can make precise the sense in which Ricci Yang-Mills solitons are weaker than Ricci solitons.  For this observation, we need the following theorem from \cite{Eber07}.

\begin{thm}\label{thm: Ricci soliton and geod flow invr} Let $N_C$ be a 2-step nilpotent Lie group of type $(p,q)$ with left-invariant metric corresponding to $C \in \mathfrak{so}(q)^p$.  The metric nilpotent group $N_C$ is both a Ricci soliton and `geodesic flow invariant' if and only if $m_1(C)=r\ Id_q$ and $m_2(C)=s\ Id_p$ for some $r,s\in \mathbb R$.
\end{thm}

\begin{cor} If $N_C$ admits a geodesic flow invariant Ricci soliton, then such a metric is also a Ricci Yang-Mills soliton.
\end{cor}

Notice that, in general, a manifold that admits a Ricci soliton will not necessarily admit a Ricci Yang-Mills soliton, as Ricci Yang-Mills solitons are only defined on manifolds that are also principal bundles.

\subsection{Ricci Yang-Mills solitons  on nilpotent Lie groups are always strong}\label{section: RYM sol of Lie type}

In this section we show that Ricci Yang-Mills solitons on nilpotent Lie groups   are  strong Ricci Yang-Mills solitons (cf. Definition \ref{defin: strong RYM soliton}).  This is true for nilpotent Lie groups of arbitrary steps, not just 2-step nilpotent.

Let $N$ be a simply-connected nilpotent Lie group with central subgroup $Z$.  The usual projection \linebreak $\pi : N \to N/Z$ is a principal bundle with structure group $G\simeq Z$.  Here $Z$ is connected and so $N/Z$ is also simply-connected.  At the identity element $e\in N$, $\pi_*:\mathfrak N \to \mathfrak N/\mathfrak Z$ has $Ker \ \pi_* = \mathfrak Z$, and  the restriction $\pi_*|_\mathcal H : \mathcal H = \mathfrak Z^\perp \to \mathfrak N / \mathfrak Z$ is a linear isometry.

\begin{lemma}  The projection $\pi$ is a Lie group homomorphism and thus
    $$\pi \circ exp_N = exp_{N/Z} \circ \pi_*. $$
Moreover, since $exp$ is a diffeomorphism we also have
    $$log_{N/Z} \circ \pi = \pi_* \circ log_N,  $$
where $log$ is the inverse of $exp$.
\end{lemma}
%

Let $\varphi \in \mathfrak{Diff}(N/Z)$ be a diffeomorphism on the base $N/Z$.  Then we define $\psi \in \mathfrak{Diff}_Z N$ by
    $$\psi  \circ exp_N ( X+ Z)= exp_N ( \pi_*|_\mathcal H ^{-1} \circ log_{N/Z} \circ \varphi \circ exp_{N/Z} \circ \pi_* (X) + Z  ),$$
where $X\in \mathcal H$ and $Z\in \mathfrak Z$.  Notice that $\psi$ is well-defined as $N$ being simply-connected implies $exp :\mathfrak N \to N$ is a diffeomorphism. By the Campbell-Baker-Hausdorff formula (cf. Theorem \ref{thm: CBH formula}) we see that $exp(Y)exp(Z) = exp(Y+Z)$ for any $Y\in \mathfrak N$ and $Z\in \mathfrak Z$, and hence  this map is a bundle automorphism.   Using the above lemma, it is straight-forward to show that $\psi$ induces $\varphi$, that is, $\pi \circ \psi = \varphi \circ \pi$.  Thus we have the following proposition.

\begin{prop}  The map $\psi$ induces $\varphi$ and hence Ricci Yang-Mills solitons on nilpotent Lie group are strong Ricci Yang-Mills solitons.
\end{prop}

\noindent
Notice that for a general principal bundle $\pi : P \to M$, one would not expect a diffeomorphism on $M$ to lift to a bundle automorphism on $P$.  The above proposition works for our nilpotent groups since the quotient is a homomorphism of Lie groups.

\section{Examples}\label{section: examples}

We now provide examples of nilpotent Lie groups that do not admit Ricci solitons but that do admit Ricci Yang-Mills solitons.  From the perspective of Geometric Invariant Theory,  that such examples exist is not a surprise.  It should be rare but not too uncommon because there should be plenty of points $C\in \mathfrak{so}(q)^p$ whose $SL_q\mathbb R\times SL_p\mathbb R$-orbit is not closed but whose $SL_q\mathbb R$-orbit is closed (cf. Corollary \ref{cor: closed SL q C orbit implies RYM soliton} and Theorem \ref{thm: Ricci soliton and geod flow invr}).

\begin{example}
\end{example}
The first example can be produced from the work of Cynthia Will.  In \cite{Wi}, Will constructs a curve of (pairwise) non-isomorphic nilpotent Lie groups which do not admit Ricci solitons.  We use the algebra corresponding to $t=1$ in her curve $\overline \mu_t$.  This algebra has structure matrices
    $$C^1 = \left[ \begin{BMAT}{cc.cc.cc}{cc.cc.cc}
        & a^2&&&&\\
        -a^2&&&&&\\
        &&&&& 1\\
        &&&&-1&\\
        &&& 1&&\\
        &&-1&&&
    \end{BMAT} \right] ,
    C^2 = \left[ \begin{BMAT}{cc.cc.cc}{cc.cc.cc}
        &&&&& a\\
        &&&&-a&\\
        &&&\ 0&&\\
        &&\ 0&&&\\
        & a&&&&\\
        -a&&&&&
    \end{BMAT} \right] ,
    C^3 = \left[ \begin{BMAT}{cc.cc.cc}{cc.cc.cc}
        &&& a&&\\
        &&-a&&&\\
        & a&&&&\\
        -a&&&&&\\
        &&&&&\ 0\\
        &&&&\ 0&
    \end{BMAT} \right]$$
The algebra presented above is isomorphic to Will's example but has different structure matrices.  The above is $g\cdot C= gCg^t$ where $g=diag\{a,a,1,1,1,1\}$ for Will's set of structure matrices $C$.

A simple computation shows that
\[m_1(C)=-2\left[ \begin{BMAT}{cc.cc.cc}{cc.cc.cc}
        -a^4-2a^2&&&&&\\
        &-a^4-2a^2&&&&\\
        &&-1-a^2&&&\\
        &&&-1-a^2&&\\
        &&&&-1-a^2&\\
        &&&&&-1-a^2
    \end{BMAT} \right]
\]

By Theorem~\ref{thm: RYM soliton vs distinguished pt}, the above will be an Ricci Yang-Mills soliton if we can show $C$ is a distinguished point; that is, if $m_1(C)\cdot C = rC$ for some $r\in \mathbb R$.  This is possible for $a^2 = \frac{-1+\sqrt 5}{2}>0$, and we have the desired result.  Notice that in this case, our soliton is generated by the trivial vector field (i.e., $D=0$ in Eqn. \ref{rym 2}), since $m_1(C)$ is a multiple of the identity.  We will provide examples below which are not generated by trivial vector fields.

\begin{example}
\end{example}
The second family of examples uses the manifolds constructed in \cite{Jab:Moduli}.  In particular, one can construct  continuous familes of algebras of type $(p,q)$ for $2\leq p \leq 6$ that are shown to not admit Ricci soliton metrics.  To do so, we must first describe a process called concatenation.

Consider  $A=(A_1, \dots , A_{p})\in \mathfrak{so}(q_1)^p$ and $B=(B_1,\dots ,B_{p})\in \mathfrak{so}(q_2)^{p}$ which are structure matrices associated to nilpotent Lie algebras $N_A$ and $N_B$ of types $(p,q_i)$, respectively.  Then we can build a new nilpotent Lie algebra $N_C$ corresponding to the structure matrix $C\in \mathfrak{so}(q)^p$, where $q=q_1+q_2$ and $$C_i=\begin{pmatrix}A_i \\ & B_i \end{pmatrix}.$$
We call this process \textit{concatenation} and denote it by $C= A+_cB$.  As $A$ and $B$ have linearly independent components, the same is true for  $C$ and hence $C$ corresponds to a nilalgebra of type $(p,q)$.  Additionally, we will abuse notation and concatenate $A\in \mathfrak{so}(q_1)^{p_1}$ and $B\in \mathfrak{so}(q_2)^{p_2}$ where $p_1 < p_2$.  This is an element of $\mathfrak{so}(q_1+q_2) ^{p_2}$ defined as
    $$ (A_1,\dots , A_{p_1}, \underbrace{0,\dots,0}_{p_2-p_1} ) +_c (B_1,\dots, B_{p_2} ).$$
We are interested in concatenating the following structure matrices.
Denote by $J$ the $2\times 2$ matrix $\begin{bmatrix} 0 & 1\\ -1 &0\end{bmatrix}$.  Define $A_1 \in \mathfrak{so}(2k)$ to be the concatenation $A_1 = J\underbrace{+_c\dots +_c}_k J$.  This is just a block diagonal matrix with all blocks being copies of $J$.  Define $B_1,B_2,\dots, B_6 \in \mathfrak{so}(4)$ as
    $$B_1= \left[ \begin{BMAT}{cc.cc}{cc.cc} 0 &1 &&\\ -1&0 &&\\ && 0&1\\ && -1&0\end{BMAT} \right], \
    B_2 = \left[ \begin{BMAT}{cc.cc}{cc.cc} && 0&\ 1 \ \\ && \ 1\ &0\\ 0&-1&&\\ -1&0&&\end{BMAT} \right], \
    B_3 = \left[ \begin{BMAT}{cc.cc}{cc.cc} && \ 1 \ &0 \\ &&0& \ 1\ \\ -1&0&&\\ 0&-1&&\end{BMAT} \right],
    $$

    $$B_4= \left[ \begin{BMAT}{cc.cc}{cc.cc} 0 &1 &&\\ -1&0 &&\\ && 0&-1\\ && 1&0\end{BMAT} \right], \
    B_5 = \left[ \begin{BMAT}{cc.cc}{cc.cc} && 0&\ 1 \ \\ && -1\ &0\\ 0&\ 1&&\\ -1&0&&\end{BMAT} \right],\
    B_6 = \left[ \begin{BMAT}{cc.cc}{cc.cc} && \ 1 \ &0 \\ &&0& -1\ \\ -1&0&&\\ 0&\ 1&&\end{BMAT} \right].
    $$
Now define $C=(C_1, \ldots, C_j)=a_1A_1+_c(b_1B_1,c_1B_2)+_c\cdots +_c(b_{n-1}B_1,c_{n-1}B_2)+_c(d_1B_1, \ldots, d_jB_j)$ for $b_i, c_i, d_i\in \mathbb{R}$.  As concatenations have such a simple presentation, it is easy to compute the value of $m_1$ at such an element.  For details see \cite{Jab:Moduli}. 

    $$ m_1(C) = -2\sum_i C_i^2 =    \left[\begin{BMAT}{cccccc}{cccccc}
        2a_1^2 Id_{2k}& &&&&\\
          &2(b_1^2+c_1^2)Id_4 &&&&\\
         && \ddots &&&\\
         &&& 2(b_{n-1}^2+c_{n-1}^2)Id_4&&\\
         &&&& 2(d_1^2 +\dots + d_j^2) Id_4   &\\
         &&&  & &
         \addpath{(0,5,.)rudrdlu}
         \addpath{(3,3,.)rldruddur}
    \end{BMAT}\right]$$
Then by Theorem~\ref{thm: RYM soliton vs distinguished pt}, an algebra of this type will admit a Ricci Yang-Mills soliton as long as
\[a_1^2=b_1^2+c_1^2=\cdots=b_{n-1}^2+c_{n-1}^2=d_1^2+\cdots+d_j^2.
\]
Thus we have a  $n-1$-parameter family of non-isomorphic algebras (by letting the $b_i$ vary) that admit Ricci Yang-Mills solitons but that do not admit Ricci solitons.

Notice that in this example, since $m_1(w)$ is a multiple of the identity, $D\equiv 0$. \\

\begin{example}
\end{example}
These examples of Ricci Yang-Mills solitons that are generated by non-trivial vector fields will be of types $(3,9), \ldots, (6,9)$.   Again, it is shown in \cite{Jab:Moduli} that these algebras do not admit Ricci soliton metrics.  Let $C$ be the concatenation

$$a_1 \left[ \J \right] +_c \lambda \left( \left[ \begin{BMAT}{c.c}{c.c}\J & \\&0 \end{BMAT}\right], \left[ \begin{BMAT}{c.c}{c.c}0&\\ &\J \end{BMAT} \right] \right) +_c (b_1 \ B_1,\dots,b_j \ B_j).$$

In this case, $m_1(C)\neq rId$, so if the manifold admits a Ricci Yang-Mills soliton, it will be nontrivial.

 Specifically, we compute $m_1(C)$ to be
 \[m_1(C)=-2\left[ \begin{BMAT}{c.ccc.c}{c.ccc.c}
       2a_1^2 Id_2&&&&\\
        &2\lambda^2&&&\\
        &&4\lambda^2&&\\
        &&&2\lambda^2&\\
        &&&&2(b_1^2+\cdots+b_j^2) Id_4
    \end{BMAT} \right].
 \]
Then a Ricci Yang-Mills soliton is admitted if $4a_1^2=6\lambda^2=4(b_1^2+\cdots+b_j^2)=r$.   Using the notation of Theorem~\ref{thm: RYM soliton vs distinguished pt}, we see that
\[B=\left[ \begin{BMAT}{c.ccc.c}{c.ccc.c}
       0&&&&\\
        &-\lambda^2&&&\\
        &&\lambda^2&&\\
        &&&-\lambda^2&\\
        &&&&0
    \end{BMAT} \right].
    \]
A simple computation confirms that $B$ is a stabilizer of $C$.  Thus we obtain a $j-1$-parameter family of non-trivial Ricci Yang-Mills soliton metrics on this algebra.  Most algebras should admit many non-isometric Ricci Yang-Mills solitons.

\appendix

\section{Lie derivatives and derivations of algebras}
Consider a Lie group $G$ with Lie algebra $\mathfrak{g}$.  Let $X$ be a vector field induced by a one-parameter family of automorphisms $\varphi_t$ of $G$.  Recall that $Lie\ Aut(G) = Aut(\mathfrak g)$  and $Lie\ Aut(\mathfrak g) = Der(\mathfrak g)$. Denote by $\psi_t$ the one-parameter group of $Aut(\mathfrak g)$ obtained as the differential of each $\varphi_t$ at the identity $e\in G$; that is, $d(\varphi_t)_e = \psi_t$ for each $t$.  Now write $\psi_t = exp(tD)$ for some $D\in Der(\mathfrak g)$.

We will show $\mathcal L_X Y = -D(Y)$ for left-invariant vector fields $Y$ on $G$, where $D$ is the derivation above.  This will then be applied to Lie derivatives of forms.  For example, for a left-invariant metric $g$ on $G$ evaluated at left-invariant vector fields $v,w$ we would have
    \begin{eqnarray*}(\mathcal L_Xg)(v,w) &=& \mathcal L_x (g(v,w)) - g(\mathcal L_Xv,w) - g(v,\mathcal L_Xw)\\
        &=& 0 + g(Dv,w) + g(v,Dw)\\
        &=& g( (D+D^t)v,w).
    \end{eqnarray*}
Similarly, we would have $(\mathcal L_x \Omega) (v,w) = \Omega(Dv,w)+\Omega(v,Dw)$, as $\Omega$ is left-invariant in our work.\\

As left-invariant vectors are global vector fields, we have $T(G) = G\times \mathfrak g$.  So given $p\in G$ and $v\in T_eG=\mathfrak g$, we have the identification $v_p = (p,v)$, where $v_p = \restrictto{\frac{d}{dt}}{0} [p\cdot exp(tv)]$ is the left-invariant vector field associated to $v\in \mathfrak g$.

Recall $(\mathcal L_X v)_p := - (\restrictto{\frac{d}{dt}}{0} \tilde \varphi_t v)_p$, where
    \begin{eqnarray*}(\tilde \varphi_t v)_p &=& \tilde \varphi_t (v_{\varphi_t^{-1}p}) = (\varphi_t)_* (\varphi_t^{-1}p,v) \\
    &=& \restrictto{\frac{d}{ds}}{0}[ \varphi_t( \varphi_t^{-1}(p)\cdot exp(sv) ) ] = \restrictto{\frac{d}{ds}}{0}[ p \cdot \varphi_t(exp(sv) ) ]\\
    &=& \restrictto{\frac{d}{ds}}{0}[ p \cdot exp(s \psi_tv) ] = (p,\psi_tv ).
    \end{eqnarray*}
Observe that if $v$ is a left-invariant vector field, so also is $\mathcal L_X v$.  Now we compute
    $$(\mathcal L_X v)_p := - (\restrictto{\frac{d}{dt}}{0} \tilde \varphi_t v)_p =(p,
    -\restrictto{\frac{d}{dt}}{0} \psi_t v) = (p, -\restrictto{\frac{d}{dt}}{0} exp(tD)v ) = (p, -Dv) $$
as desired.

\end{document}